\documentclass[letter,11pt]{article}
\usepackage{amsfonts}
\usepackage{amssymb}
\usepackage{amsmath}
\usepackage{cite}
\newtheorem{theo}{Theorem}
\newtheorem{prop}[theo]{Proposition}
\newtheorem{lema}[theo]{Lemma}
\newtheorem{coro}[theo]{Corollary}

\newcommand{\Q}{\mathbb{Q}}
\newcommand{\R}{\mathbb{R}}
\newcommand{\X}{\mathcal{X}}

\newcommand{\T}{\mathbb{T}^1}

\newcommand{\C}{\mathbb{C}}
\newcommand{\N}{\mathbb{N}}

\newcommand{\Z}{\mathbb{Z}}

\newcommand{\F}{\mathcal{F}}

\title{Hyperbolization of cocycles by isometries of the euclidean space}
\author{Mario Ponce}
	\date{}
\begin{document}
\maketitle
\begin{abstract}
We study hyperbolized versions of cohomological equations that appear with cocycles by isometries of the euclidean space. These (hyperbolized versions of) equations have a unique continuous solution. We concentrate in to know whether or not these solutions converge to a  genuine solution to the original equation, and in what sense we can use them as good approximative solutions. The main advantage of considering solutions to hyperbolized cohomological equations is that they can be easily described, since they are  global attractors of a naturally defined skew-product dynamics. We also include some technical results about twisted Birkhoff sums and exponential averaging. 
\end{abstract}
In dynamical systems, problems involving a lot of hyperbolicity, in the sense of that there is a good amount of expansiveness and contractiveness, are in general much more accessible  than problems where hyperbolicity is absent (for instance, in the dynamics of isometries, of indifferent fixed points, etc.) A naive, nevertheless fruitful, strategy for tackling  these non-hyperbolic  problems is to consider small perturbations of the maps such that the situation falls into a hyperbolic setting. In the case that this hyperbolized version of the problem presents good answers, it is natural to ask whether or not these provide us with some information about the original unperturbed problem, when the perturbation (hence the hyperbolicity) tends to disappear.  One call this strategy the {\it hyperbolization technique}. This technique has been exploited successfully in many circumstances. For instance, Yoccoz \cite{YOCC95} uses this technique for offering an extremely simple proof  that almost every indifferent quadratic polynomial is linearizable. In the context of cohomological equations,  Bousch \cite{BOUS00}  applies the hyperbolization technique in the proof of his expansive Ma\~n\'e's Lemma. Afterwards, Jenkinson \cite{JENK06}  employs the technique for finding normal forms in the setting of ergodic optimization. In a context that is directly related to our current work, Jorba and Tatjer \cite{JOTA08}  study the fractalization of invariant curves for the hyperbolized version of real one-dimensional translation cocycles over a rotation of the circle.   
\section{Framework and problems}
Let $T:X\to X$ be a homeomorphism of a compact metric space $X$, and $H$ be a metric space. A {\it cocycle by isometries  of $H$ over $T$} is a map $I:\Z\times X\to Isom(H)$ defined from $\Z\times X$ to the group of isometries of $H$ verifying the cocycle relation
\[
I(n+m, x)=I(n, T^mx)I(m, x)
\]
for every $n, m\in \Z$. The above relation allows to uniquely determine the cocycle  using the values of $I(\cdot):=I(1, \cdot)$ by the recursive relation
\[
I(n, x)=I(n-1,Tx)  I(x).
\]
Given a subgroup $\mathcal{G}\subset Isom(H)$, one of the most important problems in the study of cocycles is to know whether or not  there exists a map $B:X\to Isom(H)$  so that the conjugacy $B(Tx)\circ I(x)\circ B(x)^{-1}$ is a cocycle taking values in $\mathcal{G}$. The study of cocycles (and this problem in particular) is closely related with the study of the following induced skew-product dynamical system 
\begin{eqnarray*}
F:X\times H&\longrightarrow& X\times H \\
(x, v)&\longmapsto& (Tx, I(x)v).
\end{eqnarray*}
In this work we are concerned with the case of continuous cocycles acting on the euclidian space $\R^l$, for $l\geq 1$. Every isometry $I\in Isom(\R^l)$ can be written (in a unique way) as $Iv=\Psi v+\rho$, with $\Psi\in U(l)$, the orthogonal group  $\R^l$, and $\rho\in \R^l$. Hence, a cocycle by isometries will be defined by two continuous functions 
\[
\Psi:X\to U(l)\quad , \quad \rho:X\to \R^l
\]
so that $I(x)v=\Psi(x)v+\rho(x)$ for every $v\in \R^l$.
Let us explore the conjugacy problem proposed above in the case $\mathcal G=U(l)$. The existence of a continuous (resp. measurable) conjugacy $B:X\to Isom(\R^l)$ so that $B(Tx)\circ I(x)\circ B(x)^{-1}$  belongs to $U(l)$ for every $x\in X$ is equivalent to the existence of a continuous (resp. measurable) invariant section for the induced skew-product dynamics $F$. Namely, if there exists such $B$, then we can conjugate $F$ by $\mathcal B (x, v):=(x, B(x)v)$  in order to obtain a skew-product dynamics in the form
\[
\mathcal B \circ F \circ \mathcal B (x, v)^{-1}=(Tx, \tilde \Psi(x)v), 
\]
for some continuous (resp. measurable) function $\tilde \Psi:X\to U(l)$. The zero section $X\times \{0\}$ is invariant by $\mathcal B \circ F\circ  \mathcal B ^{-1}$ and then the section $\mathcal B^{-1}(X\times \{0\})$ is a continuous (resp. measurable) invariant section for $F$. In the other hand, if there exists $u:X\to \R^l$ a continuous (resp. measurable) invariant section for $F$, that is, $F(x, u(x))=(Tx, u(Tx))$ for every $x\in X$, then we can write $B(x)v:=v-u(x)$. A simple computation yiels 
\[
B(Tx)I(x)B(x)^{-1}0=0 \quad \textrm{for all} \quad x\in X.
\]
Hence, the isometry $B(Tx)I(x)B(x)^{-1}$ belongs to $U(l)$ and $B$ is the desired continuous (resp. measurable) conjugacy that solves the problem. The existence of a continuous  (resp. measurable) invariant section $u:X\to \R^l$ for $F$ is equivalent to the existence of a continuous (resp. measurable) solution to the {\it twisted cohomological equation}
\begin{equation}\label{tce}
u(Tx)-\Psi(x)u(x)=\rho(x)\quad \textrm{for all }x\in X.
\end{equation}
If we assume that there exists a continuous invariant section for $F$, then it is easy to see that every orbit remains  a constant distance away from  the invariant section and in particular, that every orbit is bounded. In \cite{CNP}, D. Coronel, A. Navas and M. Ponce  shown that under the additional assumption that $T$ is a minimal homeomorphism, the existence of a continuous invariant section is equivalent to the existence of a bounded orbit. In the case  that a continuous invariant section exists one says that the cocycle is a {\it continuous coboundary}. There is a big amount of results concerning the existence, regularity and rigidity of solutions to the untwisted version of the cohomological equation (that is, $\Psi\equiv id_{\R^l}$). The survey \cite{KARO} by A. Katok and E.A. Robinson is a mandatory reference in the field. However, the classical Gottschalk-Hedlund's Theorem \cite{GOHE} is the first reference whenever $T$ is minimal. In case $T$ is hyperbolic there is a whole line of research concerning the  {\it periodic orbit obstruction}, that started with the work by Li\v vsic  \cite{LIVS} and has been extended to more complex target groups by many authors (see for instance \cite{DELAWI}, \cite{KALI} and references therein contained). In the last years, a lot of attention has been payed to the resolution of the untwisted cohomological equation in the context of partially hyperbolic dynamics (mainly motivated by the works by A. Katok and A. Kononenko \cite{KAKO} and A. Wilkinson \cite{WILK}). 
\\

Another situation that also allows to get a good amount of information about the dynamics of $F$ is that of the existence of a sequence $\{u_n:X\to \R^l\}_{n\geq 0}$ of continuous sections that are {\it almost invariant}, in the sense that 
\begin{equation}\label{lrc}
\lim_{n\to \infty} |u_n(Tx)- \Psi (x)u_n(x)-\rho(x)|=0
\end{equation}
uniformly in $x\in X$.  In this case one says that $F$ is a coboundary  in {\it reduced cohomology}. 
\\

\noindent
The aim of this work is to apply a hyperbolization technique to these two situations (\ref{tce}) and (\ref{lrc}). For $\lambda\in [0,1)$ we consider the {\it $\lambda$-hyperbolized twisted cohomological equation}
\begin{equation}\label{ltce}
\lambda u_{\lambda}(Tx)-\Psi(x)u_{\lambda}(x)=\rho(x)\quad \textrm{for all }x\in X.
\end{equation}
By clearing terms and iterating (\ref{ltce}), we deduce that there exists a unique (continuous) solution:
\begin{eqnarray}\nonumber
u_{\lambda}(x)&=&-\Psi(x)^{-1}\rho(x)+\lambda\Psi(x)^{-1}u_{\lambda}(Tx)\\\nonumber
&=&-\Psi(x)^{-1}\rho(x)+\lambda\left(-\Psi(x)^{-1}\Psi(Tx)^{-1}\rho(Tx)+\lambda\Psi(x)^{-1}\Psi(Tx)^{-1}u_{\lambda}(T^2x)\right)\\\nonumber
&\vdots&\\\label{sol_ltce}
&=&-\sum_{j\geq 0} \lambda^j\Psi(x)^{-1}\cdots\Psi(T^jx)^{-1}\rho(T^jx).
\end{eqnarray}
Notice that $u_{\lambda}$ is continuous since $|\Psi(x)^{-1}\cdots\Psi(T^jx)^{-1}\rho(T^jx)|=|\rho(T^jx)|\leq \sup_{x\in X}|\rho(x)|$ for every $j$.  
In precise terms, in this work we will treat the following problems:
\\

\noindent{\bf  Problem A. } Assume that there exists a bounded solution to  (\ref{tce}). Determine  conditions in order to establish that $\lim_{\lambda\to 1^-}u_{\lambda}$ exists and equals a solution to (\ref{tce}). Determine the type of convergence  that occurs. 
\\

\noindent{\bf  Problem B. } Determine  conditions in order to establish that
\begin{equation*}
\lim_{\lambda\to 1^{-}} |u_{\lambda}(Tx)- \Psi (x)u_{\lambda}(x)-\rho(x)|=0
\end{equation*}
holds uniformly on $x\in X$. 
\\
%

The most important remark to do at this point is that in general we are not providing new information about the existence of solutions to (\ref{tce}), nor about the reducibility of the cocycle in reduced cohomology. In fact, for the second problem we will find exactly the same hypothesis that appears in the work by Bochi and Navas \cite{BONA} about reduced cohomology of cocycles. Nevertheless, we want to stress that the main advantage of working with the hyperbolized version (\ref{ltce}) of the twisted cohomological equation is that the continuous section $u_{\lambda}$ can be easily described since it is a global attractor for the  inverse of the induced skew-product dynamics:
\[
F_{\lambda}(x, v)=(Tx, \lambda^{-1}(\Psi(x)v+\rho(x))).
\]
Indeed, we can state the following
\begin{prop}
For every $\lambda \in (0, 1)$, the section $u_{\lambda}$ is invariant under $G_{\lambda}:=F_\lambda^{-1}$. Moreover, for every $(x, v)\in X\times \R^l $ one has
\[
\lim_{n\to \infty}|u_{\lambda}(T^{-n}x)- G_{\lambda}^{n}(x, v)|=0.
\]
\end{prop}
\noindent{\it Proof. } First notice that $G_{\lambda}(x, v)=(T^{-1}x, \lambda \Psi(T^{-1}x)^{-1}v-\Psi(T^{-1}x)^{-1}\rho(T^{-1}x))$. Obviusly $u_{\lambda}$ is invariant for $G_{\lambda}$ since it is invariant for $F_{\lambda}$. Hence, we can conjugate $G_{\lambda}$ by $(x, v)\mapsto (x, v-u_{\lambda}(x))$ to obtain a cocycle by linear conformal maps $(x, v)\mapsto (T^{-1}x,\lambda \Psi(T^{-1}x)^{-1}v)$. From this we can deduce that 
\[
|u_{\lambda}(T^{-n}x)- G_{\lambda}^{n}(x, v)|=\lambda^{n}|v|. \quad_{\square}
\]
\vspace{0.4cm}

This paper is organized in the following way. In sections \ref{sec_erftbs} and \ref{sec_etba} we present some technical results concerning the convergence of twisted Birkhoff sums and exponential averaging. We use these results to tackle Problem A. (in section \ref{problem1}) and Problem B. (in section \ref{problem2}).  The main results of this work are Theorem \ref{teo19} and Theorem \ref{prop1}.

\section{Ergodic results for twisted Birkhoff sums}\label{sec_erftbs}
Let $T:X\to X$ be a homeomorphism of a compact metric space $X$. Consider a continuous map $\F:X\to U(l)$ from $X$ to the topological group of linear isometries of $\R^l$. We define an extension of $T$ by
\begin{eqnarray*}
\tilde T:X\times U(l)&\longrightarrow&X\times U(l)\\
(x, A)&\longmapsto& (Tx, A\F(x)).
\end{eqnarray*} 
Given a measurable function $f:X\to \R^l$, we define the following measurable function
\begin{eqnarray*}
\tilde f:X\times U(l)&\longrightarrow& \R^l\\
(x, A)&\longmapsto& Af(x).
\end{eqnarray*}
We are interested in the Birkhoff sums of $\tilde f$ over $\tilde T$:
\begin{equation} \label{one_star}
\frac{1}{N}\sum_{j=0}^{N-1}\tilde f(\tilde T^j(x, A))
=A\left[\frac{1}{N}\sum_{j=0}^{N-1}\F(x)\F(Tx)\cdots\F(T^{j-1}x)f(T^jx)\right].
\end{equation}
The quantity above between brackets is called a {\it twisted Birkhoff sum}. The purpose of this section is to understand the behavior of these sums, and to determine various sufficient conditions in order to obtain some kind of convergence for them. Given a $\tilde T$-invariant probability measure $\tilde \mu$, the projection measure 
\[
\mu(\mathfrak{B}):=\tilde \mu (\mathfrak{B}\times U(l)), 
\]
defined for any borel set $\mathfrak{B}\subset X$,  is a $T$-invariant probability measure. Thus, if there exists some convergence result  related to the Birkhoff sums $\frac{1}{N}\sum_{j=0}^{N-1}\tilde f(\tilde T^j(x, A))$, for $\tilde \mu$ almost every pair $(x, A)\in X\times U(l)$, then the corresponding result holds for the twisted Birkhoff sums, for $\mu$ almost every point $x\in X$.  For example, we can state the following
\begin{prop}\label{p_5711}
If $\tilde T$ is uniquely ergodic and $f\in C^0(X, \R^l)$, then the following limit holds for every $x\in X$ 
\[
\lim_{N\to \infty}\frac{1}{N}\sum_{j=0}^{N-1}\F(x)\F(Tx)\cdots\F(T^{j-1}x)f(T^jx)=0
\] 
and the convergence is uniform in $x\in X$.
\end{prop} 
\noindent{\it Proof. } A limit in (\ref{one_star}) and the unique ergodicity of $\tilde T$ yield
\[
\lim_{N\to \infty} A\left[\frac{1}{N}\sum_{j=0}^{N-1}\F(x)\F(Tx)\cdots\F(T^{j-1}x)f(T^jx)\right]=\int_{X\times U(l)} \tilde fd\mu,
\] 
where $\mu$ is the unique $\tilde T$-invariant probability measure, and the limit is uniform in $x, A$. Since it holds for every $A\in U(l)$ we conclude that $\int \tilde fd\mu=0$ and the result follows.$\quad_{\square}$ 
\\

Notice that the above result also holds if we only assume that there exists a sufficiently large compact subset $K\subset U(l)$ such that $X\times K$ is $\tilde T$-invariant and $\tilde T\big |_{X\times K}$ is uniquely ergodic.
\paragraph{Invariant sections and invariant functions for $\tilde T$.} The study of the extension $\tilde T$ is closely related to the following equation 
\begin{equation}\label{two_star}
f(Tx)=\F(x)^{-1}f(x)\quad \textrm{for all}\quad x\in X, \quad f:X\to \R^l.
\end{equation}
\begin{lema}\label{lema13}
A function $f:X\to \R^l$ is a solution to (\ref{two_star}) if and only if the function $\tilde f(x, A):=Af(x)$ is $\tilde T$-invariant.
\end{lema}
\noindent{\it Proof. } One has
\begin{equation*}
\tilde f\circ \tilde T(x, A)=\tilde f(Tx, A\F(x))=A\F(x)f(Tx).
\end{equation*}
Hence 
\[
\tilde f\circ \tilde T(x, A)=\tilde f(x, A)\!\quad\! \iff\quad A\F(x)f(Tx)=Af(x)\!\quad\! \iff \quad f(Tx)=\F(x)^{-1}f(x).\quad_{\square}
\]
\begin{coro}
Assume that the unique $\tilde T$-invariant continuous functions are the constant functions (for example, this is the case when  $\tilde T$ is uniquely ergodic). The equation (\ref{two_star}) only admits the trivial solution $f(x)\equiv 0$. 
\end{coro}
\noindent{\it Proof.} Let $f$ be a solution for (\ref{two_star}). Lemma \ref{lema13} implies that $Af(x)=c$ for some $c\in \R^l$ and every $A, x$. Since $A$ ranges in $U(l)$ we conclude that $c=0$ and then $f(x)=0$ for every $x\in X. \quad_{\square}$
\\

Now, we turn to the problem of the existence of an invariant section of $\tilde T$, that is, the problem of finding a map $\X:X\to U(l)$ so that $\tilde T(x, \X(x))=(Tx, \X(Tx))$. This is equivalent to finding a solution for 
\begin{equation}\label{three_star}
\X(Tx)=\X(x)\F(x)\quad \textrm{for all} \quad x\in X, \quad \X:X\to U(l).
\end{equation}
Equations (\ref{two_star}) and (\ref{three_star}) can be related by the following
\begin{lema}\label{lema15}
If the equation (\ref{three_star}) admits a solution then equation (\ref{two_star}) admits non-trivial solutions.
\end{lema}
\noindent{\it Proof. } Let $\X:X\to U(l)$ a solution for (\ref{three_star}). Define $f(x):=\X(x)^{-1}e$, for some $e\in \R^{n}\setminus \{0\}$. We can compute
\begin{equation*}
f(Tx)=\X(Tx)^{-1}e=(\X(x)\F(x))^{-1}e=\F(x)^{-1}\X(x)^{-1}e=\F(x)^{-1}f(x),
\end{equation*}
that is, $f$ is a nontrivial solution for (\ref{two_star}).$\quad_{\square}$
\\

Let us relate the existence of $\tilde T$-invariant sections with the convergence of the twisted Birkhoff sums. Let $\X:X\to U(l)$ be a continuous $\tilde T$-invariant section. Consider the restriction $\hat T:\mathcal{G}\to \mathcal{G}$ of $\tilde T$ to the graph $\mathcal{G}:=\{(x, \X(x))\ : \ x\in X\}\subset X\times U(l)$. The invariant measures of $\hat T$ are in correspondence with the invariant measures of $T$. Let $\mu$ be a $T$-invariant measure and $\hat \mu$ the corresponding $\hat T$-invariant measure. For a continuous $f:X\to \R^l$ we define $\hat f:\mathcal{G}\to \R^l$ by $\hat f(x, \X(x))=\X(x)f(x)$ which is continuous. The classical ergodic theorem allows to say that, for $\mu$ almost every point, the following (pointwise) limit exists
\[
g(x):=\lim_{N\to \infty}\frac{1}{N}\sum_{j=0}^{N-1}\hat f\circ \hat T^j(x, \X(x))=\X(x)\left[\lim_{N\to \infty }\frac{1}{N}\sum_{j=0}^{N-1}\F(x)\cdots \F(T^{j-1}x)f(T^jx)\right].
\]
Furthermore, the function $g$ is $T$-invariant. Since $\X$ is a solution to (\ref{three_star}), the invariance of $g$ gives 
\begin{eqnarray*}
\X(Tx)\left[\lim_{N\to \infty }\frac{1}{N}\sum_{j=0}^{N-1}\F(Tx)\cdots \F(T^{j}x)f(T^{j+1}x)\right]&=&\X(x)\left[\lim_{N\to \infty }\frac{1}{N}\sum_{j=0}^{N-1}\F(x)\cdots \F(T^{j-1}x)f(T^jx)\right]\\
\F(x)\left[\lim_{N\to \infty }\frac{1}{N}\sum_{j=0}^{N-1}\F(Tx)\cdots \F(T^{j}x)f(T^{j+1}x)\right]&=&\lim_{N\to \infty }\frac{1}{N}\sum_{j=0}^{N-1}\F(x)\cdots \F(T^{j-1}x)f(T^jx)
\end{eqnarray*}
for $\mu$ almost every $x\in X$. Hence, the function 
\[
x\longmapsto \lim_{N\to \infty }\frac{1}{N}\sum_{j=0}^{N-1}\F(x)\cdots \F(T^{j-1}x)f(T^jx)
\]
is a solution to the equation (\ref{two_star}), for $\mu$ almost every $x\in X$. As an application, we can show the following 
\begin{prop}\label{p_5713}
Let $\X:X\to U(l)$ be a continuous  $\tilde T$-invariant section and $f:X\to \R^l$ a continuous function. If $T$ is uniquely ergodic then the following limit exists 
\[
h(x):=\lim_{N\to \infty }\frac{1}{N}\sum_{j=0}^{N-1}\F(x)\cdots \F(T^{j-1}x)f(T^jx), 
\]
and it is uniform in $x\in X$. Moreover, for $c=\int_X \X(x)f(x)d\mu$ one has $h(x)=\X(x)^{-1}c$, that is a continuous solution to the equation (\ref{two_star}).$\quad_{\square}$
\end{prop}
\noindent{\bf Example I.} Take $X=\T$, $T(\theta)=\theta+\alpha$, for some $\alpha \in \R\setminus \Q$.  We consider $l=2,\  \R^2\sim \C$ and for some fixed $\beta\in \R$ we put $\F(\theta)\equiv e^{i\beta}$ (that is, $\F(\theta)$ is the constant rotation by angle $\beta$). The restriction of $\tilde T$ to the compact $\tilde T$-invariant set $\T\times \mathbb{S}^1$ acts as $\tilde T(\theta, e^{it})=(\theta+\alpha, e^{i(t+\beta)})$. We have two cases:
\begin{itemize}
\item If $(\alpha, \beta)$ are rationally independent, then $\tilde T\big |_{\T\times \mathbb{S}^1}$ is uniquely ergodic.  Hence, Proposition \ref{p_5711} applies and for every $f\in C^0(\T, \C)$ one has 
\begin{equation*}
\lim_{N\to \infty} \frac{1}{N}\sum_{j=0}^{N-1}e^{ij\beta}f(\theta+j\alpha)=0
\end{equation*}
uniformly in $\theta$.
\item If $(\alpha, \beta)$ are not rationally independent, then (perhaps a finite covering of) $\tilde T\big |_{\T\times \mathbb{S}^1}$  has a continuous invariant section. Hence, Proposition \ref{p_5713} applies and for every $f\in C^0(\T, \C)$ the following limit
\begin{equation*}
\lim_{N\to \infty} \frac{1}{N}\sum_{j=0}^{N-1}e^{ij\beta}f(\theta+j\alpha)
\end{equation*}
exists, uniformly on $\theta\in \T. \quad_{\triangle}$
\end{itemize}
\noindent{\bf Remark. } Assume that $T$ is uniquely ergodic, and that $f$ is continuous. The limit $h(x)$ (which by Proposition \ref{p_5712} is defined for almost every $x$ with respect to the unique $T$-invariant measure) is not necessarily uniform. Indeed, one can modifies the Example I. taking a continuous non-constant $\mathcal{F}(\theta)=e^{i\beta(\theta)}$ such that the cohomological equation $\psi(\theta+\alpha)-\psi(\theta)=\beta (\theta)$ has a measurable a.e. defined solution $\psi$ but no continuous solution (as in the Furstenberg example). Then a calculation using Birkhoff theorem shows that $h(\theta)=ce^{-i\psi(\theta)}$ for a.e. $\theta$, where $c=\int e^{i\psi(\theta)}f(\theta)d\theta$. If we choose a continuous $f$ such that $c\neq 0$, then the function $h$ is not continuous, and in particular the limit that defines it cannot be uniform. Hence, the hypotheses of Propositions \ref{p_5711} and \ref{p_5713} are necessary. 
\paragraph{Invariant measures for $\tilde T$.} Since $U(l)$ is a topological group, there exists a probability measure $\mathcal{L}$ in $U(l)$ that is invariant under right multiplication in the group (the {\it Haar} measure). Given a $T$-invariant probability measure $\mu$ the {\it product probability measure} 
\[
 \mu_{\mathcal L} (\mathfrak B\times \mathfrak A):= \mu(\mathfrak B)\mathcal{L}(\mathfrak A)
\]
defined for every borel set $\mathfrak B\times \mathfrak A\subset X\times U(l)$, is $\tilde T$-invariant (see for instance \S  3.6 in \cite{BOGO}). Notice that in particular this yields:
\begin{lema}
If $\tilde T$ is uniquely ergodic then $T$ is uniquely ergodic.$\quad_{\square}$
\end{lema}
The classical Birkhoff ergodic theorem allows then to obtain 
\begin{prop}\label{p_5712}
 Let $\mu$ be a $T$-invariant measure and $f:X\to \R^l$ be a measurable function so that $\int |f(x)|d\mu <\infty$ (we say $f\in L^1_\mu(X, \R^l)$). Then for $\mu$ almost every $x\in X$ the following (pointwise) limit exists
\[
h(x):=\lim_{N\to \infty }\frac{1}{N}\sum_{j=0}^{N-1}\F(x)\cdots \F(T^{j-1}x)f(T^jx).
\]
Moreover, $h$ belongs to $L^1_\mu(X, \R^l)$, and the convergence occurs in $L^1_\mu(X, \R^l)$. Finally, $h$ is a solution to the equation (\ref{two_star}), for $\mu$ almost every $x\in X.$
\end{prop}
\noindent{\it Proof. } Since $f\in L^1_\mu(X, \R^l)$, the function $\tilde f(x, A)=Af(x)$ belongs to $L^1_{\mu_{\mathcal{L}}}(X\times U(l), \R^l)$. The Birkhoff ergodic theorem applied to $\tilde f$ and to the $\tilde T$-invariant probability measure $\mu_{\mathcal{L}}$ allows to conclude.$\quad_{\square}$
\paragraph{The Mean Ergodic Theorem by von Neumann.} Let $H$ be a real Hilbert space, $\mathfrak U:H\to H$ be a weak linear contraction, that is, $| \mathfrak Uf|_H\leq |f|_H$ for every $f\in H$. Let $M=\{f\in H \ |\ \mathfrak Uf=f\}$ be the closed subspace of the invariant elements by $\mathfrak U$ and $P:H\to M$ the projection. The Mean Ergodic Theorem by von Neumann asserts that for every $f\in H$
\[
\frac{1}{N}\sum_{j=0}^{N-1} \mathfrak U^jf
\]
converges (in the norm $|\cdot|_H$), to the projection $Pf$. As a direct application we can state the following:
\begin{prop}\label{571_14}
Let $\mu$ be a $T$-invariant measure and $f:X\to \R^l$ be a function in $L^2_\mu(X, \R^l)$. Then there exists $h\in L^2_\mu(X, \R^l)$ so that 
\[
\lim_{N\to \infty} \left\|h(x)-\frac{1}{N}\sum_{j=0}^{N-1}\F(x)\cdots \F(T^{j-1}x)f(T^jx)\right\|_{L^2_\mu(X, \R^l)}=0.
\]
Moreover, the function $h$ verifies 
\[
h(Tx)=\F(x)^{-1}h(x)\quad \mu-a.e.
\]
\end{prop}
\noindent{\it Proof. } We pick $H:=L^2_\mu(X, \R^l)$, and the linear map $\mathfrak U:L^2_\mu(X, \R^l)\to L^2_\mu(X, \R^l)$ defined as
\[
\mathfrak Uf(x):=\F(x)f(Tx)\quad \textrm{for} \quad f\in L^2_\mu(X, \R^l), x\in X.
\]
Since $\F(x)$ is a linear isometry and $\mu$ is $T$-invariant we have
\[
\| \mathfrak Uf\|_{L^2_\mu(X, \R^l)}^2=\int_X |\F(x)f(Tx)|^2d\mu=\int_X |f(Tx)|^2d\mu=\int_X |f(x)|^2d\mu=\|f\|^2_{L^2_\mu(X, \R^l)}.
\]
Hence we can apply the von Neumann Mean Ergodic Theorem. The result easily follows by noticing that for every $j$ one has
\[
\mathfrak U^jf(x)=\F(x)\cdots\F(T^{j-1}x)f(T^jx). \quad_{\square}
\]
\section{Exponential twisted Birkhoff averages}\label{sec_etba}
\paragraph{Summability methods and a theorem by Frobenius.} Consider a sequence of vectors $\{z_j\}_{j\geq 0}$ in a real  normed  vector space $V$. We say that the sequence $\{z_j\}$ is {\it convergent by the Ces\`aro $C(1)$ method} to the value $\tilde z\in V$ if
\[
\lim_{n\to \infty}\frac{1}{n}\sum_{j=0}^{n-1}z_j=\tilde z.
\]
Assume that the formal power series $(1-\lambda)\sum_{j=0}^{\infty}\lambda^jz_j$ has convergence radius at least $1$. Then we say that the sequence $\{z_j\}$ is {\it convergent by the Abel method} to the value $\tilde z\in V$ if
\[
\lim_{\lambda\to 1^{-}}(1-\lambda)\sum_{j=0}^{\infty}\lambda^jz_j=\tilde z.
\]
\begin{theo}[Frobenius, see \cite{PETE66}]
If the sequence $\{z_j\}$ is convergent by the Ces\`aro $C(1)$ method to the value $\tilde z$ then it is convergent by the Abel method to the same value $\tilde z$.
\end{theo}
\noindent{\bf Remark. } Notice that if $\{z_n\}$ converges by the Ces\`aro $C(1)$ method, then $\frac{z_n}{n}\to 0$, and in particular the power series $\sum \lambda^jz_j$ has convergence radius at least $1$.\\ 
We are going to provide a  proof of this theorem in the next paragraphs. Before this, and as a manner  of converse theorem we mention the following
\begin{theo}[Tauberian Theorem, see \cite{PETE66}, \cite{HARD92}]
If $\{z_j\}$ is bounded and convergent by the Abel method then it is convergent by the Ces\`aro $C(1)$ method.$\quad_{\blacksquare}$
\end{theo}
In order to prove the Frobenius Theorem, for every $n\in \N$ and any sequence $\{a_j\}_{j\geq 0} \subset V$ we define 
\[
s_{-1}=0\quad , \quad s_n=a_0+a_1+\dots+a_{n}.
\]
For $\lambda\in (0, 1)$, and assuming that $\sum \lambda^jz_j$ has convergence radius at least $1$,  we compute
\begin{eqnarray}\nonumber
(1-\lambda)\sum_{n\geq 0}\lambda^n s_n&=&\sum_{n\geq 0}\lambda^n s_n-\sum_{n\geq 0}\lambda^{n+1} s_n =\sum_{n\geq 0}\lambda^n (s_n-s_{n-1})\\ \label{estrella1}
&=& \sum_{n\geq 0}\lambda^n a_n.
\end{eqnarray}
Notice that (\ref{estrella1}) shows that $\sum \lambda^ns_n$ has convergence radius at least $1$. 
For every $n\in \N$ define $\sigma_n=\frac{1}{n}(s_0+s_1+\dots+s_{n-1})$.
\begin{lema}[Frobenius]\label{lemafrobenius}
\[
\textrm{If }\quad  \lim_{n\to \infty}\sigma_n=A\quad  \textrm{ then }\quad  
\lim_{\lambda\to 1^{-1}}\sum_{n\geq 0}\lambda^na_n=A.
\]
\end{lema}
\noindent{\it Proof. } Using (\ref{estrella1}) twice we obtain
\begin{eqnarray*}
\sum_{n\geq 0}\lambda^n a_n&=&(1-\lambda)\sum_{n\geq 0}\lambda^n s_n=(1-\lambda)^2\sum_{n\geq 0}\lambda^n(s_0+s_1+\dots+s_{n})\\
&=&(1-\lambda)^2\sum_{n\geq 0}\lambda^n(n+1)\sigma_{n+1}.
\end{eqnarray*}
Using that $(1-\lambda)^2=\left(\sum_{n\geq 0}(n+1)\lambda^n\right)^{-1}$ we get
\begin{eqnarray}\nonumber
\sum_{n\geq 0}\lambda^n a_n&=& \frac{\sum_{n\geq 0}\lambda^n(n+1)\sigma_{n+1}}{\sum_{n\geq 0}(n+1)\lambda^n}\\\label{55_1}
&=&A+\frac{\sum_{n\geq 0}\lambda^n(n+1)(\sigma_{n+1}-A)}{\sum_{n\geq 0}(n+1)\lambda^n}\\\nonumber
&\stackrel{\lambda\to1^{-}}{\longrightarrow}&A. \quad_{\square}
\end{eqnarray}
Now we can easily obtain the proof of the Frobenius Theorem by putting $a_0=z_0$, $a_j=z_j-z_{j-1}$. In this way we have
\[
s_n=z_n\quad , \quad {\sigma_n}=\frac{1}{n}\sum_{j=0}^{n-1}z_j \quad , \quad (1-\lambda)\sum_{n\geq 0 }\lambda^n s_n= (1-\lambda)\sum_{n\geq 0 }\lambda^n z_n
\] 
and the theorem follows from Lemma \ref{lemafrobenius} and equality (\ref{estrella1}).$\quad_{\square}$
\\

The following is a functional version of the Frobenius Theorem. Let $(X, \mu)$ be a probability space, $\{z_j:X\to V\}_{j\geq 0}$ a sequence of measurable $V$-valued functions belonging to $L^1_\mu(X, V)$, that is $\int_X |z_j(x)|_Vd\mu(x)<\infty$. Define $\sigma_n(x)=\frac{1}{n}(z_0(x)+z_1(x)+\dots+z_{n-1}(x))$. Using equality (\ref{55_1}) one easily obtains 
\begin{lema}\label{lema_571}
\noindent
\begin{itemize}
\item[$(a)$ ] If there exists $A\in L^1_\mu(X, V)$ such that $\sigma_n(x)\to A(x)$ pointwise for $\mu$-a.e. $x\in X$ and in $L^1_\mu(X, V)$, whenever $n\to \infty$, then the following convergence holds pointwise for $\mu$-a.e. $x\in X$ and in $L^1_\mu(X, V)$ 
\[
(1-\lambda)\sum_{j\geq 0}\lambda^jz_j(x)\to A(x)\quad \textrm{when}\quad \lambda\to 1^{-}.
\]  
\item[$(b)$ ] If there exists $A\in L^1_\mu(X, V)$ such that $\sigma_n(x)\to A(x)$ uniformly in $x$ whenever $n\to \infty$ and $\{z_j:X\to V\}_{j\geq 0}$ is uniformly bounded, then the following convergence holds uniformly in $x\in X$
\[
(1-\lambda)\sum_{j\geq 0}\lambda^jz_j(x)\to A(x)\quad \textrm{when}\quad \lambda\to 1^{-}.\quad_{\square}
\] 

\end{itemize}
\end{lema}
\paragraph{Exponential averaging and some related ergodic results.} In classical ergodic theory, given a measurable dynamics $T:X\to X$ and a measurable observable $f:X\to \R$, we are concerned with studying the asymptotic behavior of the Birkhoff averages $\frac{1}{n}\sum_{0\leq j< n} f\circ T^j(x)$. Notice that this average  is a (total mass $1$) weighted sum of the values taken by the observable over a finite piece of the orbit $x, Tx, \dots, T^{n-1}x$. Thus, the classical ergodic theorems deal with the limits of these averages whenever $n$ tends to $\infty$. Notice that in this situation, the fact of letting $n\to \infty$ can be interpreted as a way to considering the values of $f$ over large pieces of the orbit, with equal weights. \\

It is then natural to extend the notion of asymptotic averages to other averaging methods. Here we propose to study the following {\it exponential averaging}. In contrast with the classical one, in this method we are going to consider the values taken by $f$ over the entire orbit $\{x, Tx, T^2x, \dots \}$. The weights we are going to consider are given by the powers of a given real number $\lambda \in (0,1)$. More precisely, we define the $\lambda$-average of $f$ over $T$ as
\begin{equation}\label{ultimomani}
S_{\lambda}(f, x)=\frac{1}{\sum_{j\geq 0}\lambda^j}\sum_{j\geq 0}\lambda^jf(T^jx).
\end{equation} 
For a fixed $\lambda$, the $\lambda$-average $S_{\lambda}(f, x)$ assigns large weights to the first iterates of $x$, and very small weights to  large iterates of $x$.  Since we are interested in the asymptotic behavior of $f$ when taking values over the entire orbit, it is natural to study the behavior of $S_{\lambda}$ as $\lambda$ grows to $1$. Our previous discussion about summability methods allows us to conclude that the classical Birkhoff averaging method covers this exponential averaging method. More precisely, we can state the next  two Propositions, that follow directly from the classical Birkhoff Ergodic Theorem  together with  Lemma \ref{lema_571}. 
\begin{prop}[Birkhoff Ergodic Theorem for exponential averaging] Let $(X, \mu)$ be a probability space, $T:X\to X$ a measure-preserving transformation and $f:X\to \R^l$ a measurable function in $L^1_\mu(X, \R^l)$, that is, $\int_X|f(x)|d\mu <\infty$. Then
\begin{enumerate}
\item $\lim_{\lambda\to 1^-}S_{\lambda}(f, x)=h(x)$ exists $\mu$-a.e.;
\item $h(Tx)=h(x)$ $\mu$-a.e.;
\item $h\in L^1_\mu(X, \R^l)$;
\item If  $\mathcal{A}\subset X$ is a $T$-invariant Borel set, then $\int_\mathcal{A}fd\mu=\int_\mathcal{A}hd\mu$;
\item $S_{\lambda}(f, \cdot)\to h$ in $L^1_\mu(X, \R^l)$, when $\lambda\to 1^-$.$\quad_{\square}$
\end{enumerate}
\end{prop}
\begin{prop} Let $T:X\to X$ be a uniquely ergodic homeomorphisms of a compact metric space and $f:X\to \R^l$ be a continuous function. Then 
\[
\lim_{\lambda\to 1^{-}}S_{\lambda}(f, x)=\int_{X}f d\mu, 
\]
where $\mu$ is the unique $T$-invariant probability measure, and the convergence is uniform in $x\in X$.$\quad_{\square}$
\end{prop}
\paragraph{Twisted exponential averages.} The Birkhoff averaging method for twisted Birkhoff sums that we treat in section \ref{sec_erftbs} can be now easily extended to exponential twisted averages 
\[
(1-\lambda)\sum_{j\geq 0}\lambda^j\F(x)\cdots \F(T^{j-1}x)f(T^jx).
\]
Namely, the previous results allows us to state the following
\begin{prop}\label{prop_comodin}
Propositions \ref{p_5711},  \ref{p_5713}, \ref{p_5712} and  \ref{571_14} also hold  whenever we replace the Birkhoff sum limits 
\[
\lim_{N\to \infty}\frac{1}{N}\sum_{j=0}^{N-1}\F(x)\cdots \F(T^{j-1}x)f(T^jx)
\]
    by exponential averages limits      
    \[
    \lim_{\lambda\to 1^{-}}(1-\lambda)\sum_{j\geq 0}\lambda^j\F(x)\cdots \F(T^{j-1}x)f(T^jx).
    \]
     Furthermore, the resulting limit functions are the same for both averaging methods.$\quad_\square$
\end{prop} 
\section{Solutions to hyperbolized equations converge to a genuine solution}\label{problem1}
Assume there exists a bounded function $u^*:X\to \R^l$ that, in some sence (eg. continuous, pointwise, a.e.),  is a solution to the twisted cohomological equation (\ref{tce}).  We are interested in (the existence of) the convergence 
\begin{equation}\label{convergence_571}
\lim_{\lambda \to 1^-}u_{\lambda, }
\end{equation}
where $u_{\lambda}$ is the unique continuous solution to the hyperbolized cohomological equation (\ref{ltce}). Notice that whenever the convergence in (\ref{convergence_571}) occurs in a given sense (eg. pointwise, uniform, a.e., $L^1$, $L^2$, etc.), then necessarily the limit object must satisfy the cohomological equation (\ref{tce}) in the corresponding way. Remember that $
u_{\lambda}(x)=-\sum_{j\geq 0}\lambda^j\Psi(x)^{-1}\cdots\Psi(T^jx)^{-1}\rho(T^jx).$ Using that $u^*$ is a bounded solution to (\ref{tce}), a direct computation yields
\begin{eqnarray}\nonumber
u_{\lambda}(x)&=&-\sum_{j\geq 0}\lambda^j\Psi(x)^{-1}\cdots\Psi(T^jx)^{-1}[u^*(T^{j+1}x)-\Psi(T^jx)u^*(T^jx)]\\  \label{18_5_dos_estrellas}
&=& u^*(x)-(1-\lambda)\sum_{j\geq 0}\lambda^{j}\Psi(x)^{-1}\cdots\Psi(T^jx)^{-1}u^*(T^{j+1}x).
\end{eqnarray}
Then, the convergence in (\ref{convergence_571}) is equivalent to the convergence of the exponential twisted average in (\ref{18_5_dos_estrellas}). Notice that the equality (\ref{18_5_dos_estrellas}) above holds in the same sense in which $u^*$ is a solution to (\ref{tce}). 
\\

Every bounded solution $u^{\star}$ to the twisted cohomological equation (\ref{tce}) has the form $u^{\star}(x)=u^*(x)+v(x)$, for some bounded function $v:X\to \R^n$ that is a solution to the homogeneous equation
\begin{equation}\label{homogenea}
v(Tx)-\Psi(x)v(x)=0.
\end{equation}
\begin{lema}
If there exists a bounded solution $u^*$ to the twisted cohomological equation (\ref{tce}), such that 
\[
(1-\lambda)\sum_{j\geq 0}\lambda^{j}\Psi(x)^{-1}\cdots\Psi(T^jx)^{-1}u^*(T^{j+1}x)
\]
 converges (when $\lambda\to 1^{-}$),  then for every bounded solution $u^{\star}$ to the twisted cohomological equation (\ref{tce}) 
 \[
 (1-\lambda)\sum_{j\geq 0}\lambda^{j}\Psi(x)^{-1}\cdots\Psi(T^jx)^{-1}u^{\star}(T^{j+1}x)
 \]
 converges (when $\lambda\to 1^{-}$). Furthermore, whenever $\mathcal{S}_{\lambda}(u^*, x)$ converges, the limit function is a solution to the homogeneous equation (\ref{homogenea}).  
\end{lema}
\noindent{\it Proof. } The operator  
\[
\mathcal{S}_{\lambda}(f, x):=(1-\lambda)\sum_{j\geq 0}\lambda^{j}\Psi(x)^{-1}\cdots
\Psi(T^jx)^{-1}f(T^{j+1}x)
\] 
is linear in $f$ and one can easily check that $\mathcal{S}_{\lambda}(v, x)=v(x)$ whenever $v$ is a solution to the homogenous equation (\ref{homogenea}). Hence we have
\[
\mathcal{S}_{\lambda}(u^{\star}, x)=\mathcal{S}_{\lambda}(u^{*}+v, x)=\mathcal{S}_{\lambda}(u^{*}, x)+v(x),
\]
and the result follows.$\quad_\square$
\\

Let $u:X\to \R^l$ a bounded function. We say that $u$ has {\it (uniform, pointwise, a.e., etc.)-exponential zero mean} if the following  
\[
\lim_{\lambda\to 1^{-}}\mathcal{S}_{\lambda}(u, x)=0,
\]
holds for every $x\in X$ (uniformly, pointwise, a.e., etc.).
\begin{lema}
The solutions $u_{\lambda}$ converge (uniformly, pointwise, a.e., etc.) when $\lambda\to 1^{-}$, if and only if there exists a (uniform, pointwise, a.e., etc.)-exponential zero mean bounded solution $u^{*}$ to the equation (\ref{tce}). In this case, $u_{\lambda}$ converges to a (uniform, pointwise, a.e., etc.)-exponential  zero mean solution to the equation (\ref{tce}).
\end{lema}
\noindent{\it Proof. } Assume that $u_{\lambda}$ converge to $u$. Then we can use $u$ for computing $u_{\lambda}$ as in (\ref{18_5_dos_estrellas}):
\[
u_{\lambda}(x)=u(x)-\mathcal{S}_{\lambda}(u, x).
\]
Taking $\lambda\to 1^{-}$ we deduce that $u$ is a exponential zero mean solution to (\ref{tce}). Conversely, assume that $u$ is a exponential zero mean solution to (\ref{tce}). We can use $u$ for computing $u_{\lambda}$ as above, and to deduce that $u_{\lambda}$ converges to $u$.$\quad_{\square}$  
\\

\noindent{\bf Example II.} Let $T:X\to X$ be a uniquely ergodic homeomorphism, $\Psi(x):=id_{\R^l}$ and $\rho:X\to \R^l$ be a continuous function. Assume that there exists a continuous solution $u^{\star}$ to the cohomological equation
\[
u(Tx)-u(x)=\rho(x).
\]
For instance, this is the case whenever the induced skew-product dynamics has a bounded orbit and $T$ is minimal (this is the Gottschalk-Hedlund's Theorem). It is known that all the solutions to this equation have the form $u^{\star}+c$ with $c\in \R^l$. Define the solution $u^*=u^{\star}-\int_X u^{\star}d\mu$, where $\mu$ is the unique $T$-invariant probability measure.   Proposition \ref{prop_comodin} asserts that $\mathcal{S}_{\lambda}(u^{*}, x)\to \int_{X}u^{*}d\mu=0$ uniformly in $x$, when $\lambda \to 1^{-}$. Hence, the solutions  to the hyperbolized cohomological equation converge uniformly to the  unique zero mean solution to the cohomological equation, that is
\[
u_{\lambda}\stackrel{\lambda\to{1^-}}{\longrightarrow} u^{*} \ \textrm{uniformly, and } \int_X{u^*}d\mu=0.\quad_{\triangle}
\] 

In the sequel, we will apply the Proposition \ref{prop_comodin} to more general situations. In order to fit with the framework of sections \ref{sec_erftbs} and \ref{sec_etba} we define $\F(x):=\Psi(x)^{-1}$ and   $f(x):=u^*(Tx)$, for some $u^*:X\to \R^l$ that is a bounded solution to the cohomological equation (\ref{tce}). Recall that we have already defined the extension 
\[
\tilde T(x, A):=(Tx, A\Psi(x)^{-1}), \quad\textrm{for }\quad (x, A)\in X\times U(l).
\]
The exponential averaging versions of the Propositions  \ref{p_5711} and  \ref{p_5713} contained in the Proposition \ref{prop_comodin} give directly
\begin{theo}\label{teo19}
supose that the equation (\ref{tce}) has a continuous solution and that $T$ is uniquely ergodic. 
\begin{itemize}
\item[$(a)$] If the extension $\tilde T$ is uniquely ergodic then there exists a continuous map $u^*:X\to \R^l$ such that
\[
\lim_{\lambda\to 1^-}u_{\lambda}(x)=u^*(x), 
\]
uniformly in $x\in X$. Moreover, $u^*$ is a uniform-exponential zero mean continuous solution to the twisted cohomological equation (\ref{tce}).
\item[$(b)$] If there exists a continuous $\tilde T$-invariant section  then there exists a continuous map $u^*:X\to \R^l$ such that
\[
\lim_{\lambda\to 1^-}u_{\lambda}(x)=u^*(x), 
\]
uniformly in $x\in X$. Moreover, $u^*$ is a uniform-exponential zero mean continuous solution to the twisted cohomological equation (\ref{tce}).$\quad_{\square}$
\end{itemize}
\end{theo}
\noindent{\bf Example III. } The following class of cocycles was initially considered in \cite{CONAPO2}. For example, in that work the authors construct topologically transitive cocycles inside this class. Let $\alpha\in \R\setminus \Q$ and $\beta\in \T$. Let $\rho:\T\to \C$ be a continuous function. We construct the cocycle by isometries of $\C$ (a {\it cylindrical vortex})
\begin{eqnarray*}
F_{\alpha, \beta, \rho}:\T\times \C&\longrightarrow& \T\times \C\\
(\theta, z)&\longmapsto&(\theta+\alpha, e^{i\beta}z+\rho(\theta)).
\end{eqnarray*}
Assume that there exists a continuous solution $u:X\to \R^l$ to the twisted cohomological equation 
\begin{equation}\label{pintura}
u(\theta+\alpha)-e^{i\beta}u(\theta)=\rho(\theta), 
\end{equation}
 Namely, this is the case whenever there exists a bounded orbit for $F_{\alpha, \beta, \rho}$ (see \cite{CNP}). The discussion in the Example I. and Theorem \ref{teo19} allows to conclude that $u_{\lambda}$ converges uniformly to a continuous solution to the twisted cohomological equation (\ref{pintura}). Also notice that in the case that the pair $(\alpha, \beta)$ is rationally independent, this solution is unique.$\quad_{\triangle}$ 
\\

\noindent
The exponential averaging versions of the Propositions  \ref{p_5712} and  \ref{571_14} contained in the Proposition \ref{prop_comodin} give directly
\begin{theo}Let $\mu$ be a $T$-invariant probability measure and suppose  that there exists  a bounded solution  to the equation (\ref{tce}),  for $\mu$-a.e. $x\in X$. Then we have:
\begin{itemize} 

\item[$(a)$]   There exists $u^*:X\to \R^l\in L_\mu^1(X, \R^l)$ that is a $\mu$-a.e. solution to the twisted cohomological equation (\ref{tce}) and such that when $\lambda\to 1^-$
\[
u_{\lambda}(x)\to u^*(x) \quad \textrm{pointwise for}\quad \mu-a.e. \quad x\in X.
\] 
Moreover, the convergence above occurs  in $L^1_{\mu}(X, \R^l)$ and  the function $u^*$ has $\mu$-a.e. exponential zero mean.
\item[$(b)$] There exists $u^*:X\to \R^l\in L_\mu^2(X, \R^l)$ that is a $\mu$-a.e. solution to the twisted cohomological equation (\ref{tce}) and such that when $\lambda\to 1^-$
\[
u_{\lambda}\to u^*\quad\textrm{in}\quad L^2_{\mu}(X, \R^l). \quad_{\square}
\] 
\end{itemize}
\end{theo}
\section{Solutions to hyperbolized equations almost solve the cohomological equation}\label{problem2}
 Let $v:X\to \R^l$ be a section. One can define the {\it displacement} of $v$ by the cocycle $I$ as
\[
Disp(v, I)=\sup_{x\in X}|v(Tx)- I(x)v(x)|.
\]
In terms of the displacement, a section $v$ is invariant by the induced skew-product $F$ if and only if $Disp(v, I)=0$. Moreover, a sequence of sections $\{v_n:X\to \R^l\}_{\N}$ is almost invariant if and only if $Disp(v_n, I)\to 0$.  
\\

 In a more general framework, Bochi and Navas \cite{BONA}  characterize the cocycles by isometries that  have a sequence of continuous almost invariant sections. They made it by considering the {\it maximal drift} of the cocycle $D(T, I)$ defined as
\[
D(T, I):=\lim_{n\to \infty} \frac{1}{n}\sup_{x\in X}dist(v_0, I(n, x)v_0)=\lim_{n\to \infty} \frac{1}{n}\sup_{x\in X}|v_0- I(n, x)v_0|
\] 
for some $v_0\in \R^n$. This limit always exists by a sub-additive argument and it is independent of the choice of $v_0$.  Whenever $D(T, I)=0$ one says that the cocycle $I$ verifies the {\it zero drift condition.}
\\

\noindent{\bf Example IV. } Assume that $F$ has a bounded invariant section $v:X\to \R^l$. Pick $v_0\in \R^l$ and compute
\begin{eqnarray*}
|v_0- I(n, x)v_0| &\leq&  |v_0- v(x)|+|v(x)- I(n, x)v(x)|+|I(n, x)v(x)- I(n, x)v_0|\\
&=&|v_0- v(x)|+|v(x)- v(T^nx)|+|v(x)-v_0|
\end{eqnarray*} 
that is uniformly bounded. Hence $D(T, I)=0$.$\quad_{\triangle}$
\\

\noindent{\bf Example V. } Assume that $F$ has a sequence of bounded almost invariant sections $v_k:X\to \R^l$. supose that $D(T, I)=d>0$. Pick $k>0$ so that $|v_k(Tx)- I(x)v_k(x)|<d/2$ for every $x\in X$. Let us compute the maximal drift by picking $v_0=v_k(x)$ for some $ x\in X$:
\begin{eqnarray*}
|v_k( x)- I(n, x)v_k( x)|&\leq& |v_k(x)- v_k(T^nx)|+|v_k(T^nx)- I(n, x)v_k(x)|\\
&\leq& |v_k(x)- v_k(T^nx)|+\\
&&+\sum_{j=0}^{n-1}|I^j(T^{n-j}x)v_k(T^{n-j}x)- I({j+1},T^{n-(j+1)}x)v_k(T^{n-(j+1)}x)|.\\
\end{eqnarray*}
But, $I(j+1, T^{n-(j+1)}x)=I(j, T^{n-j}x)I(T^{n-(j+1)}x)$ and then 
\begin{eqnarray*}
|v_k(x)- I(n, x)v_k(x)|&\leq& |v_k(x)- v_k(T^nx)|+\\
&&+\sum_{j=0}^{n-1}|v_k(T\circ T^{n-(j+1)}x)- I(T^{n-(j+1)}x)v_k(T^{n-(j+1)}x)|\\
&\leq&  |v_k(x)- v_k(T^nx)|+n\frac{d}{2}, 
\end{eqnarray*}
which is impossible for large $n$. Then we conclude that $D(T, I)=0$.$\quad_{\triangle}$
\\


In fact, the next theorem of Bochi and Navas \cite{BONA} is the equivalence between the zero drift condition and the existence of a sequence of continuous almost invariant sections for the cocycle.
\begin{theo}[Bochi and Navas, see Theorem A in \cite{BONA}]
Assume that $H$ is a Busemann space. Let  $I:X\to Isom(H)$ be a cocycle over $T:X\to X$. If $D(T, I)=0$ then the cocycle is a coboundary in reduced cohomology.$\quad_{\blacksquare}$
\end{theo}
The proof of this theorem relies on the construction of a certain type of {\it barycenter} in $H$.  In this work we propose to show a version of this theorem in the case of $H=\R^l$ by making use of the hyperbolization technique. 
\\

In order to obtain precise estimations, let us compute some iterates of the cocycle:
\begin{eqnarray*}
I(x)v&=&\Psi(x)v+\rho(x)\\
I(2, x)v&=&\Psi(Tx)\Psi(x)v+\Psi(Tx)\rho(x)+\rho(Tx)\\
&\vdots&\\
I(n, x)v&=&\left(\prod_{j=0}^{n-1}\Psi(T^jx)\right)v+\sum_{j=0}^{n-1}\left(\prod_{k=j+1}^{n-1}\Psi(T^kx)\right)\rho(T^jx).
\end{eqnarray*}
In particular we have 
\[
\Psi(x)^{-1}\cdots \Psi(T^{n-1}x)^{-1}I(n, x)v=v+\sum_{j=0}^{n-1}\Psi(x)^{-1}\cdots\Psi(T^jx)^{-1}\rho(T^jx).
\]
We can compute the maximal drift by picking $v_0=0$:
\begin{eqnarray*}
D(T, I)&=&\lim_{n\to \infty}\frac{1}{n}\sup_{x\in X}| I(n, x)0|\\
&=&\lim_{n\to \infty}\frac{1}{n}\sup_{x\in X}|\Psi(x)^{-1}\cdots \Psi(T^{n-1}x)^{-1}0- \Psi(x)^{-1}\cdots \Psi(T^{n-1}x)^{-1}I(n, x)0|\\
&=&\lim_{n\to \infty}\frac{1}{n}\sup_{x\in X}\left | \sum_{j=0}^{n-1}\Psi(x)^{-1}\cdots\Psi(T^jx)^{-1}\rho(T^jx) \right |.
\end{eqnarray*}
The reader should notice that the above limit fits into the framework of the twisted Birkhoff sums treated in section \ref{sec_erftbs}.     In this way the zero drift condition reads
\begin{lema}\label{lema_promedios_ergodicos}
The maximal displacement $D(T, I)$ is zero if and only if 
\[
\lim_{n\to \infty}\frac{1}{n} \sum_{j=0}^{n-1}\Psi(x)^{-1}\cdots\Psi(T^jx)^{-1}\rho(T^jx) =0
\]
uniformly in $x\in X$.$\quad_{\square}$
\end{lema}
Let us relate the solutions to the hyperbolized cohomological equations (\ref{ltce}) with the zero drift condition. Remember that the solution $u_{\lambda}$ can be explicitly  written as
\[
u_{\lambda}(x)=-\sum_{j\geq 0}\lambda^j\Psi(x)^{-1}\cdots\Psi(T^jx)^{-1}\rho(T^jx).
\]
The section $u_\lambda$ fails to be invariant for $F$ since $I(x)u_{\lambda}(x)=\lambda u_{\lambda}(Tx)\neq u_{\lambda}(Tx)$. In fact, the {displacement} of the section $u_{\lambda}$ is
\begin{eqnarray*}
Disp(u_{\lambda}, I)&=&\sup_{x\in X}|u_{\lambda}(Tx)-I(x)u_{\lambda}(x)|\\ 
&=&\sup_{x\in X}|(1-\lambda)u_{\lambda}(Tx)|=\sup_{x\in X}\left|\frac{\sum_{j\geq 0}\lambda^j\Psi(x)^{-1}\cdots\Psi(T^jx)^{-1}\rho(T^jx)}{\sum_{j\geq 0}\lambda^j}\right|.
\end{eqnarray*}
Hence, we can state the following lemma which should be compared to Lemma \ref{lema_promedios_ergodicos}.
\begin{lema}\label{lema_17mayo_18}
The displacement $Disp(u_{\lambda}, I)$ goes to $0$ when $\lambda \to 1^-$ if and only if
\[
\lim_{\lambda\to 1^- }\frac{\sum_{j\geq 0}\lambda^j\Psi(x)^{-1}\cdots\Psi(T^jx)^{-1}\rho(T^jx)}{\sum_{j\geq 0}\lambda^j}=0.
\]
uniformly in $x\in X$.$\quad_{\square}$
\end{lema}
Notice that (like in the expression of Lemma \ref{lema_promedios_ergodicos}), the above expression  is a weighted sum of the terms $\Psi(x)^{-1}\cdots\Psi(T^jx)^{-1}\rho(T^jx)$ with total mass equal to $1$. As seen in the Example V. above, the condition $D(T, I)=0$ is a necessary condition in order to get almost invariant bounded  sections. The next result states that this is a sufficient  condition in order to  the hyperbolization technique  provides almost invariant continuous sections.
\begin{theo}\label{prop1}
If $D(T, I)=0$ then 
\[
\lim_{\lambda\to 1^-} |u_{\lambda}(Tx)-\Psi(x)u_{\lambda}(x)-\rho(x)|=0
\]
 uniformly in $x\in X$. 
\end{theo}
\noindent{\it Proof.} We need to prove that $Disp(u_{\lambda}, I)$ goes to $0$ when $\lambda \to 1^-$. Putting 
\[
z_j(x):= \Psi(x)^{-1}\cdots\Psi(T^jx)^{-1}\rho(T^jx), 
\]
 Lemma \ref{lema_17mayo_18} says that we need to prove that 
\[
\lim_{\lambda \to 1^{-}}(1-\lambda)\sum_{j\geq 0}\lambda^jz_j(x)=0
\]
uniformly in $x$.  Lemma \ref{lema_571} asserts that in order to get this, it is enough to prove 
\[
\lim_{n\to \infty}\frac{1}{n}\sum_{j=0}^{n-1}z_j(x)=0
\]
uniformly in $x$. Indeed, this is the content of Lemma \ref{lema_promedios_ergodicos}.$\quad_{\square}$
\\

\noindent{\bf Example VI. } Let $(\alpha , \beta)\in \T\times \T$ be a rationally independent pair of angles and  $\rho:\T\to \C$ be a continuous function. Consider the cocycle $F_{\alpha, \beta, \rho}$ as defined at the Example III.   In order to determine the maximal drift we need to compute
\[
\lim_{n\to \infty}\frac{1}{n}\sum_{j=0}^{n-1}e^{-ij\beta}\rho(\theta+j\alpha).
\]
We have already seen at the Example I. that this limit vanishes uniformly in $\theta$, and hence the cocycle $F_{\alpha, \beta, \rho}$ verifies the zero drift condition.$\quad_{\triangle}$
\\

\begin{small}

\noindent{\bf Acknowledgments.} The preparation of this manuscript was partially funded by  the Fondecyt Grant 11090003 and the project {\it Center of Dynamical Systems and Related Fields} ACT1103. This work was motivated during fruitful discussions with Daniel Coronel, Godofredo Iommi, Jan Kiwi, Renaud Leplaideur and  Andr\'es Navas.  The author acknowledges the comments and remarks made by the Referee in order to improve the presentation of this article.

\end{small}

\begin{footnotesize}

\vspace{0.25cm}

\noindent{Mario Ponce}

\noindent{Facultad de Matem\'aticas, Universidad Cat\'olica de Chile}

\noindent{Casilla 306, Santiago 22, Chile}

\noindent{E-mail: mponcea@mat.puc.cl}

\end{footnotesize}

 \end{document}